# NO NONLOCALITY. NO FRACTIONAL DERIVATIVE.

**Vasily E. Tarasov**,
Skobeltsyn Institute of Nuclear Physics,
Lomonosov Moscow State University, Moscow 119991, Russia
E-mail: tarasov@theory.sinp.msu.ru

**Abstract:** The paper discusses the characteristic properties of fractional derivatives of non-integer order. It is known that derivatives of integer orders are determined by properties of differentiable functions only in an infinitely small neighborhood of the considered point. Therefore differential equation, which is considered for this point and contains a finite number of integer-order derivatives, cannot describe nonlocality in space and time. This allows us to propose a principle of nonlocality for fractional derivatives. We state that if the differential equation with fractional derivative can be presented as a differential equation with a finite number of integer-order derivatives, then this fractional derivative cannot be considered as a derivative of non-integer order. This means that all results obtained for this type of fractional derivatives can be derived by using differential operators with integer orders. To illustrate the application of the nonlocality principle, we prove that the conformable fractional derivative, the M-fractional derivative, the alternative fractional derivative, the local fractional derivative and the Caputo-Fabrizio fractional derivatives with exponential kernels cannot be considered as fractional derivatives of non-integer orders.

**Keywords**: fractional derivative, nonlocality, memory, conformable fractional derivative, local fractional derivative, Caputo-Fabrizio fractional derivative

**MSC:** 26A33  Fractional derivatives and integrals

**MSC:** 34A08  Fractional differential equations

**Introduction**

In modern mathematics, the fractional derivatives of non-integer orders are well-known. There are many different types of fractional derivatives (for example, see [1,2,3,4]). These fractional derivatives have been suggested by such famous mathematicians as Riemann, Liouville, Letnikov, Sonin, Weyl, Riesz, Erdelyi, Kober and other. The fractional derivatives of non-integer orders are actively used in the natural sciences to describe the processes and systems with spatial and temporal nonlocality (the nonlocality in time is usually called memory). In general, the fractional derivatives of non-integer orders form a special class of integro-differential operators. Fractional derivatives of non-integer order have a set of non-standard properties (for example, see [1,2,3,4] and [5,6,7,8,9,]). The question arises: Is it possible to formulate the main characteristic properties of fractional derivatives? As the main characteristic property, it has been proposed a violation of the standard Leibniz rule [5]. We proved that all fractional derivatives, which satisfy the standard Leibniz rule, should have the integer order, i.e. fractional derivatives of non-integer orders cannot satisfy the standard Leibniz rule [5]. However, this property is necessary, but not sufficient. There are many integro-differential operators, for which the standard Leibniz rule is violated, but they are not derivatives of non-integer order. In recent years, various types of operators have been proposed that are attempted to be classified as fractional derivatives. Despite the fact that the authors call these operators as fractional derivatives, it is important to understand whether these operators can be considered as derivatives of non-integer order, or they are simply re-designation of well-known differential operators of integer orders. Obviously, if we call the operator $D^\alpha := \alpha \cdot \frac{d}{dt}$ as a fractional derivative of non-integer order α, then this does not mean that it is such, since this operator is differential operator of first order. However, in more complex cases, it is not so obvious to understand whether the proposed operator is fractional or not. In this connection the question arises about the formulation of principles that allow us to identify fractional derivatives of non-integer orders. In this paper we propose principles for identifying fractional derivatives of non-integer order based on the following facts.

At first, it is known that the derivatives of integer orders are determined by properties of differentiable functions only in infinitely small neighborhood of the considered point. As a result, the processes or systems, which are described by this differential equation with finite number of integer-order derivatives, are considered in the infinitesimally small region and it cannot be considered as nonlocal. This means that the nonlocal effects with respect cannot be described by this equation.

Secondly, the nonlocality by time (dynamic memory) means a dependence of output (endogenous) variable at the present time on the changes of input (exogenous) variable on finite (or infinite) time interval of the past [10,11]. The memory can be considered as a property of processes, which characterizes the dependence of this process at a given time on the states in the past [10]. This means that behavior of process with memory is determined not only by the state of the process $\{t, X(t)\}$ at a given moment of time t or an infinitely small neighborhood of this point, but also by the states $\{\tau, X(\tau)\}$ at finite (or infinite) time interval (for example, $\tau \in [t_0, t]$ or $\tau \in (-\infty, t)$). To describe processes with memory, we should take into account that the output variable Y(t) and its integer-order derivatives $Y^{(k)}(\tau)$ ($k = 1,2, \ldots, m < \infty$) at time t depends on the changes of the input variable X(τ) and its integer-order derivatives $X^{(k)}(\tau)$ ($k = 1,2, \ldots, n < \infty$) for $\tau < t$ on a finite (or infinite) time interval [10,11].

Thirdly, for applied mathematics, physics and economics, the most important property of fractional derivatives of non-integer orders is that they are a new tool for describing processes and systems with spatial nonlocality and dynamic memory (nonlocality in time). Mathematically the property of nonlocality is expressed in the fact that fractional derivatives of non-integer orders cannot be represented as a finite set of derivatives of integer orders. Therefore, it is important to have a mathematical formulation of nonlocality property, which would allow us to identify fractional derivatives of non-integer order. As a basis of the formulation of nonlocality principle, we can consider the fact that differential equations, which contain only derivatives of integer orders, cannot describe nonlocal processes and systems.

To explain the proposed approach, one can consider a differential equation, which contains a fractional derivative of input (exogenous) or output (endogenous) variables with respect to time. Let us assume that this differential equation with fractional derivative can be presented in the form of a differential equation with a finite number of derivatives of integer-orders. The following question arises. Is it possible to regard this fractional derivative as a fractional derivative of non-integer order in this case? To formulate an answer to this question we can use a concept of equivalence of differential equations with respect to solutions. We can call two differential equations the equivalent on some wide function space if solutions of these equations coincide on this space.

## 2. Principles of nonlocality for fractional derivatives

Let us formulate the principles of nonlocality of fractional derivatives. The operator, which we will check whether it is a fractional derivative, will be called a tested fractional derivative. We call a tested equation such a differential equation that has the form

$$F\left((D_{t_0}^\alpha X)(t); X^{(k)}(t); Y^{(k)}(t); Y^{(k)}(t_0); t; t_0\right) = 0 \qquad (1)$$

which contains a tested fractional derivative $(D_{t_0}^\alpha X)(t)$ of the function X(t), the integer order derivatives $X^{(k)}(t)$ ($k = 0,1, \ldots, n_F < \infty$) and $Y^{(k)}(t)$ ($k = 0,1, \ldots, m_F < \infty$), the values $Y^{(k)}(t_0)$ ($k = 1,2, \ldots, m < \infty$), the variable t and an initial value $t_0$.

**Principle of nonlocality of fractional derivative:**
*If for the tested fractional derivative $D_{t_0}^\alpha$ there exists tested differential equation (1) that can be represented as a linear differential equation of the form*

$$I(X^{(k)}(t); Y^{(k)}(t); t) = 0 \qquad (2)$$

*which contains only integer-order derivatives $X^{(k)}(t)$ ($k = 0,1, \ldots, n_I < \infty$), $Y^{(k)}(t)$ ($k = 0,1, \ldots, m_I < \infty$) and the variable t, then the tested fractional derivative $(D_{t_0+}^\alpha X)(t)$ cannot be considered as a fractional derivative of non-integer order.*

Equation (2) assumes that all integer-order derivatives $X^{(k)}(t)$ and $Y^{(k)}(t)$ must exist. By virtue of this, the equivalence of (1) and (2) can be considered in the function space $C^n(\Omega)$, where $\Omega = [t_0, t_1]$ is a certain region of the real axis and $n = \max\{n_I, m_I\}$. We also can use a more narrow space of functions. For example, we can consider analytic functions in an interval $(t_0, t_1)$ of real axis, i.e., for function, which is can be represented as a power series in this interval.

We note that differential equation (2) with derivatives of integer order must be independent of the initial time $t_0$ of the initial values of the variables $X(t_0)$, $Y(t_0)$ and its derivatives $X^{(k)}(t_0)$, $Y^{(k)}(t_0)$ ($k = 0,1, \ldots, n < \infty$) at $t_0$.

Since differential equation (2) contains a finite number of derivatives of the integer order, then this equation, and consequently the tested fractional derivative, cannot describe processes and systems with nonlocality and memory.

For simplification, we can formulate the linear principle of nonlocality. We will call a linear tested equation such a differential equation that has the form

$$\lambda \cdot (D_{t_0}^\alpha X)(t) + \sum_{k=0}^{n_F} \mu_k(t) \cdot X^{(k)}(t) + \sum_{k=0}^{m_F} \left(v_k(t) \cdot Y^{(k)}(t) + \eta_k(t, t_0) \cdot Y^{(k)}(t_0)\right) = 0, \qquad (3)$$

which contains the tested operator $(D_{t_0}^\alpha X)(t)$, the integer order derivatives $X^{(k)}(t)$ ($k = 0,1, \ldots, n_F < \infty$) and $Y^{(k)}(t)$ ($k = 0,1, \ldots, m_F < \infty$), the values $Y^{(k)}(t_0)$ ($k = 1,2, \ldots, m < \infty$), the variable t, the initial value $t_0$ and $\lambda \neq 0$.

**Linear principle of nonlocality of fractional derivative:**
*If for the tested fractional derivative $D_{t_0}^\alpha$ there exists linear tested differential equation (3) that can be represented in the form of the linear differential equation*

$$\sum_{k=0}^{n_I} a_k(t) \cdot X^{(k)}(t) + \sum_{k=0}^{m_I} b_k(t) \cdot Y^{(k)}(t) + c(t) = 0, \qquad (4)$$

*which contains a finite number of integer-order derivatives, then the tested fractional derivative $(D_{t_0}^\alpha X)(t)$ cannot be considered as fractional derivative of non-integer order. Here the functions $X(t)$ and $Y(t)$ can be assumed to be analytic functions in an interval $(t_0, t_1)$ of real axis (i.e., for function, which is can be represented as a power series in this interval) or from the function space $C^n([t_0, t_1])$, where $n = max\{n_I, m_I\}$.*

The linear nonlocality principle actually means that it is impossible to represent a linear differential equation with fractional derivatives of non-integer order in the form of an equation with polynomial differential operators $P_n(D)$, where $n < \infty$, $D = d/dt$ and coefficients can be non-constant.

In the simplest case, this principle states that if there is a linear differential equation

$$\lambda \cdot (D_{t_0}^\alpha X)(t) + \nu(t, t_0) \cdot Y(t) + \eta(t, t_0) Y(t_0) = 0, \qquad (5)$$

that can be represented in the form the differential equation

$$a_1(t) \cdot \frac{dX(t)}{dt} + a_0(t) \cdot X(t) + b_0(t) \cdot Y(t) = 0, \qquad (6)$$

then the tested fractional derivative $D_{t_0}^\alpha$ cannot be considered as a fractional derivative of non-integer order.

### 3. Examples of application of nonlocality principles.

Let us give examples of application of nonlocality principle for fractional derivatives. We prove that the conformable, local and Caputo-Fabrizio fractional derivatives cannot be considered as fractional derivatives of non-integer order. Differential equations containing these operators can be described by differential equations with a finite number of derivatives of integer order. This leads to the fact that results, which are obtained by using these operators, can be derived by using the differential operators with integer orders.

**Example 1: The conformable fractional derivative**
Recently, the conformable fractional derivative $T_\alpha$ has been proposed by Khalil, Horani, Yousef and Sababheh [12] in 2014 and then were considered in different. In Theorem 2.2 [12, p.67] the standard Leibniz rule was proved for the conformable fractional derivative. This means that this derivative is a derivative of integer (first) order [5]. In Definition 2.1 of [12, p.66] the conformable fractional derivative of the non-integer order $\alpha \in (0,1)$ is defined by the equation

$$(T_\alpha X)(t) = \lim_{\varepsilon \to 0} \frac{X(t + \varepsilon \cdot t^{1-\alpha}) - X(t)}{\varepsilon}, \qquad (7)$$

where the function is defined for $t > 0$.

It is easy to verify that the differential equation

$$\lambda \cdot (T_\alpha X)(t) + \mu(t) \cdot X(t) + \nu(t) \cdot Y(t) = 0 \qquad (8)$$

with the conformable fractional derivative $T_\alpha$, where $\lambda \neq 0$ and $t > 0$, is the equivalent to the differential equation of first order of the form

$$a_1(t) \cdot \frac{dX(t)}{dt} + a_0(t) \cdot X(t) + b_0(t) \cdot Y(t) = 0, \tag{9}$$

where $a_1(t) = t^{1-\alpha} \cdot \lambda$, $a_0(t) = \mu(t)$, $b_0(t) = \nu(t)$. This statement is based on Theorem 2.2 of [12, p.67] that gives the relationship between the conformable derivative and the first derivative in the form

$$(T_\alpha X)(t) = t^{1-\alpha} \frac{dX(t)}{dt}. \tag{10}$$

Expression (10) means that the conformable fractional derivative $T_\alpha$ is a differential operator of integer (first) order with variable coefficient ($a(t) = t^{1-\alpha}$). This actually means that all results obtained for the conformable derivatives can be derived by using differential operators with integer orders.

As a result, differential equations with conformable fractional derivatives can be represented as differential equations of integer order for the space of differentiable functions. The applicability of the conformable derivatives for non-differentiable continuous functions is outside the scope of the principle of nonlocality, since this principle considers equivalence to equations with derivatives of integer orders, which assume differentiability of the function. Therefore we can state that the conformable fractional derivatives do not give anything new in the spaces of differentiable functions and are not fractional derivatives of non-integer order.

**Example 2: The fractional derivative of Katugampola**

In 2014, Katugampola [13] proposed a fractional derivative, called an alternative fractional derivative. The Katugampola (alternative) fractional derivative is defined (see Definition 2.1 in [13, p. 2] and [14]), by the equation

$$(D^\alpha X)(t) = \lim_{\varepsilon \to 0} \frac{X(t \cdot \exp(\varepsilon \cdot t^{-\alpha})) - X(t)}{\varepsilon}, \tag{11}$$

where $t > 0$. For operator (11), the standard Leibniz rule is satisfied (see Theorem 2.3 of [13, p.3]) and hence the alternative fractional derivative is a derivative of first order [5]. The relationship between operator (11) and the first derivative has the form (10) (see Theorem 2.3 of [13, p.3] and equations 2-4 in [14, p. 063502-2]). As a result, the fractional derivative (11) is a differential operator of integer (first) order with variable coefficient $a(t) = t^{1-\alpha}$. As a result, differential equations with the Katugampola (alternative) fractional derivatives can be represented as differential equations of integer order for the space of differentiable functions.

**Example 3: The M-fractional derivative**

The M-fractional fractional derivative is defined (see Definition 2.1 in [15, p. 85]), by the equation

$$(D^\alpha X)(t) = \lim_{\varepsilon \to 0} \frac{X(t \cdot E_{\beta/m}(\varepsilon \cdot t^{-\alpha})) - X(t)}{\varepsilon}, \tag{12}$$

where the function $E_{\beta/m}(z)$ is given by the equation

$$E_{\beta/m}(z) = \sum_{k=0}^{m} \frac{z^k}{\Gamma(\beta \cdot k + 1)}, \tag{13}$$

where $\beta > 0$. For operator (12), the standard Leibniz rule is satisfied (see Theorem 2.2 of [15, p.87]) and this fractional derivative is a derivative of first order [5].. In the space of differentiable functions, the relationship between operator (12) and the first derivative is given in Theorem 2.2 of [15, p.87] in the form

$$(T_\alpha X)(t) = \frac{t^{1-\alpha}}{\Gamma(\beta + 1)} \frac{dX(t)}{dt}. \tag{14}$$

As a result, the M-fractional derivative (12) is also a differential operator of integer (first) order with variable coefficient $a(t) = t^{1-\alpha}/\Gamma(\beta + 1)$. As a result, differential equations with M-fractional derivatives can be represented as differential equations of integer order for the space of differentiable functions.

**Example 4: The fractional derivative of Kolwankar and Gangal**

About 20 years ago, the local fractional derivative has been proposed by Kolwankar and Gangal in [16,17,18]. In Equation 26 of [17] and in Definition 2.2 of [18] the local fractional derivative $D_{KG}^\alpha$ is defined by the equation

$$(D_{KG}^\alpha X)(a) = \lim_{t \to a} D_{RL;a+}^\alpha \left( X(t) - \sum_{k=0}^{n-1} \frac{X^{(k)}(a)}{k!} (t - a)^k \right), \tag{25}$$

where $D_{RL;a+}^\alpha$ is the Riemann-Liouville fractional derivative and $\alpha \in (n - 1, n)$ [4, p.70-71].

In Proposition 2.1 of [19, p.730-731], [18, p.7] the standard Leibniz rule was proved for the local fractional derivative. This means that operator (15) is a derivative of integer (first) order on the space of differentiable functions [5].

It is easy to verify that the differential equation

$$\lambda \cdot (D_{KG}^\alpha X)(t) + \mu(t) \cdot X(t) + \nu(t) \cdot Y(t) = 0 \tag{36}$$

with the local fractional derivative $(D_{KG}^\alpha X)$ of the order $\alpha \in (n - 1, n)$, where $\lambda \neq 0$, is the equivalent to the differential equation of first order in the form

$$a_n(t) \cdot \frac{d^n X(t)}{dt^n} + a_0(t) \cdot X(t) + b_0(t) \cdot Y(t) = 0, \tag{47}$$

where $a_n(t) = \lambda$, $a_0(t) = \mu(t)$, $b_0(t) = \nu(t)$. This statement is based on Proposition of [20, p.197-198]. In paper [20] we proved that local fractional derivatives of differentiable functions are integer-order derivative or zero operator. If the function is n-differentiable at interval $(t_0, t_1)$, then the local fractional derivative of the order $n - 1 <$ alpha $< n$ at $t = t_0$ of this function at all points $t_0$ of $(t_0, t_1)$ is equal to zero or it is the standard derivative of integer order n. This means that the local fractional derivative is a differential operator of integer order.

As a result, differential equations with the local fractional derivatives can be represented as the differential equations of integer order for the space of differentiable functions. The problems of applying local derivatives for non-differentiable continuous functions [21] are outside the scope of the proposed principle of nonlocality, since this principle considers equivalence with standard differential equations, and hence, assumes differentiability of the function. At least in spaces of differentiable functions, local fractional derivatives do not give anything new and are not fractional derivatives of non-integer order.

**Example 5: The Caputo-Fabrizio fractional derivative**

Recently, the fractional derivatives with exponential kernels have been proposed by Caputo and Fabrizio [22] in 2015 and then were considered in various works. The Caputo-Fabrizio fractional derivative $D_{CF}^{(\alpha)}$ of the non-integer order $\alpha \in (0,1)$ is defined [22, p.74] by the equation

$$\left(D_{CF}^{(\alpha)}X\right)(t) = \frac{m(\alpha)}{1-\alpha} \cdot \int_{t_0}^{t} \exp\left\{-\frac{\alpha}{1-\alpha} \cdot (t-\tau)\right\} \cdot X^{(1)}(\tau) d\tau, \tag{58}$$

where $m(\alpha)$ is a normalization function $m(0) = m(1) = 1$. For $n > 1$, the Caputo-Fabrizio fractional derivative of the order $\alpha + n \in (n, n+1)$ is defined [22, p.76] by the expression

$$\left(D_{CF}^{(\alpha+n)}X\right)(t) = \left(D_{CF}^{(\alpha)}X^{(n)}\right)(t), \tag{69}$$

where $\alpha \in (0,1)$.

Let us consider the linear differential equation

$$\lambda \cdot \left(D_{CF}^{(\alpha)}X\right)(t) + v(t) \cdot Y(t) + \eta(t,t_0) \cdot Y(t_0) = 0, \tag{20}$$

where $\alpha \in (0,1)$, $\eta(t,t_0) = \exp\{a(\alpha) \cdot (t-t_0)\}$, $\lambda \neq 0$ and $v(t) = 1$. We can prove that differential equation (16) with the Caputo-Fabrizio fractional derivative $\left(D_{CF}^{(\alpha)}X\right)$ of the order $\alpha \in (0,1)$, where $\lambda \neq 0$, is the equivalent to the differential equation of first order

$$a_1(t) \cdot \frac{dX(t)}{dt} + b_1(t) \cdot \frac{dY(t)}{dt} + b_0(t) \cdot Y(t) = 0, \tag{27}$$

where $a_1(t) = \lambda \cdot \frac{m(\alpha)}{1-\alpha}$, $a_0(t) = a(\alpha)$, $b_1(t) = v(t)$. The proof of this equivalence is based on the well-known theorem of linear differential equations (for example, see Theorem 1.1 (Equations 1.8 and 1.9) of [23]). This theorem described solutions of the differential equation

$$\frac{dY(t)}{dt} = A(\alpha) \cdot Y(t) + F(t), \tag{22}$$

where $A(\alpha)$ and $B(\alpha)$ the constants that depend only on the parameter $\alpha \in (0,1)$. The integration of equation (22) gives the expression

$$Y(t) = \exp\{A(\alpha) \cdot (t-t_0)\} \cdot Y(t_0) + \int_{t_0}^{t} \exp\{A(\alpha) \cdot (t-\tau)\} \cdot F(\tau) d\tau. \tag{23}$$

Using the constants

$$A(\alpha) = -\frac{\alpha}{1-\alpha}, \tag{24}$$

$$B(\alpha) = \lambda \cdot \frac{m(\alpha)}{1-\alpha}, \tag{25}$$

the function $F(t) = B(\alpha) \cdot X^{(1)}(t)$, and the definition of the Caputo-Fabrizio fractional derivative (18), expression (23) can be written in the form

$$Y(t) = \eta(t,t_0) \cdot Y(t_0) + \lambda \cdot \left(D_{CF}^{(\alpha)}X\right)(t), \tag{26}$$

where $\eta(t,t_0) = \exp\{A(\alpha) \cdot (t-t_0)\}$.

Substitution of the expression (26) or/and (23) into equation (22) gives an identity. As a result, we can consider equations (22) and (26) as equivalent. Since differential equation (22) contains only the derivative of first order, equation (26) cannot describe processes and systems with nonlocality and memory. As a result, the Caputo-Fabrizio operator (18) cannot be considered as a fractional derivative of non-integer order.

Note that using the variable

$$Z(t) = Y(t) - B(\alpha) \cdot X(t), \tag{27}$$

and equation (22), equation (26) can be represented in the form

$$\frac{dZ(t)}{dt} = A(\alpha) \cdot Z(t) + C(\alpha) \cdot X(t), \tag{28}$$

where $C(\alpha) = A(\alpha) \cdot B(\alpha)$. This form of the equation allows you to interpret the action of Caputo-Fabrizio fractional derivative, as a superposition (parallel action) of the standard accelerator and the multiplier without memory [11].

For the Caputo-Fabrizio fractional derivative $D_{CF}^{(\alpha)}$ of the order $\alpha \in (1,2)$, an representation of fractional differential equation in the form of differential equation of second-order can be derived from the results in [24]. For example, fractional differential equation (1) of [24] with the Caputo-Fabrizio derivative of order $\alpha \in (1,2)$ can be represented [24,p.258] in the form of differential equation of second order (see equation 6 of [24]). The second example is given for fractional differential equation 11 of [24, p.260].

Note that in [25], it was proved that the Caputo-Fabrizio fractional derivative implements an integer order high-pass filter and that this operator is neither fractional, nor a derivative. To proof this statement, the equation $Y(t) = \left(D_{CF}^{(\alpha)} X\right)(t)$ is analyzed.

## 4. Conclusion

As a result, in Examples 1-5 we proved that the conformable fractional derivative, the alternative and M-fractional fractional derivatives, the local fractional derivative of Kolwankar and Gangal, the Caputo-Fabrizio fractional derivatives with exponential kernels cannot be considered as fractional derivatives of non-integer orders and all results obtained for this type of operators can be derived by using differential operators of integer orders. This means that all results obtained for these operators can be derived by using the differential operators with integer orders. Therefore, the proposed operators do not give anything new at least in spaces of differentiable functions except for the change of notations.

Note that the Riemann-Liouville, Caputo, Hadamard, Marchaud fractional derivatives of non-integer order have an important property that allows us to express these fractional derivatives through an infinite series of derivatives of integer orders (for example, see Lemma 15.3 of [1, p.278] and [26,27]). For example, for the space of analytic functions in an interval $(t_0, t_1)$ of real axis (i.e., for function, which is can be represented as a power series in this interval) the Riemann-Liouville fractional derivative of the order $\alpha > 0$ can be represented (see Lemma 15.3 of [1, p.278]) as an infinite series of the integer-order derivatives in the form

$$\left(D_{RL;t_0+}^{\alpha} X\right)(t) = \sum_{k=0}^{\infty} a_k(t, t_0, \alpha) \cdot X^{(k)}(t), \qquad (29)$$

where the coefficients $a_k(t, t_0, \alpha)$ are defined as

$$a_k(t, t_0, \alpha) = \binom{\alpha}{k} \cdot \frac{(t - t_0)^k}{\Gamma(n - \alpha + 1)}, \qquad (30)$$

and $\binom{\alpha}{k}$ is the generalized binomial coefficient (see equation 1.48 of [1]). This fact means that the linear differential equation (3) with fractional derivatives of non-integer order cannot be represented in the form the differential equation of finite integer order for analytic functions in the interval $(t_0, t_1)$ in general. The fact that an operator cannot be specified as a finite sum of derivatives of an integer order in general is an important sign of the nonlocality of the operator.

Because of this, the locality of an operator is synonymous with the fact that equations with this operator can be described by differential equations of integer order. In this case, nonlocality is the characteristic property of a fractional differential operator of non-integer order.

It is obvious that not all non-local operators can be considered fractional differential operators, since at least there are fractional integrals of non-integer order [1,2,3,4].